\documentclass[a4paper, 12pt]{article}
\usepackage{amssymb,amsmath,latexsym,color,graphicx}
\usepackage{enumerate}
\usepackage{color}
\usepackage[francais]{babel}

\newtheorem{dref}{D\'efinition}[section]
\newtheorem{theo}[dref]{Th\'eor\`eme}
\newtheorem{prop}[dref]{Proposition}
\newtheorem{remark}[dref]{Remarque}

\begin{document}
\date
{}
\title{ $\mathcal{PT}$ sym\'{e}trie et puits de potentiel}

\author{Naima Boussekkine \\
{\small Universit\'{e} de Mostaganem }\\
{\small Facult\'{e} des science exactes et informatique }\\
{\small 27000-Mostaganem, Alg\'erie}\\
{\small nboussekkine@yahoo.fr} \and Nawal Mecherout \\
{\small Universit\'{e} de Mostaganem }\\
{\small Facult\'{e} des science exactes et informatique }\\
{\small 27000-Mostaganem, Alg\'erie}\\
{\small mecheroutnawel@yahoo.fr } }
\date{}
\title{ PT sym\'{e}trie et puits de potentiel}
\maketitle

\begin{abstract}
Dans ce travail, on consid\'{e}re des perturbations $\mathcal{PT}$-sym\'{e}%
triques d'un op\'{e}rateur de Schr\"{o}dinger semi-classique auto-adjoint
sur la droite r\'{e}el dans le cas d'un puits de potentiel simple. On
suppose que le potentiel soit analytique et on montre que les valeurs
propres restent r\'{e}elles sous la perturbation.
\end{abstract}

\tableofcontents

\section{Introduction}

\label{intr} \setcounter{equation}{0} Les op\'{e}rateurs $\mathcal{PT}$-sym%
\'{e}triques ont a \'{e}t\'{e} propos\'{e}es comme une alternative aux op%
\'{e}rateurs auto-adjoints en physique quantique. La r\'{e}alit\'{e} du
spectre est alors importante du point de vue de la physique (voir, \cite{Be}%
, \cite{BeBeMa02} , \cite{BeMa10}, \cite{BeBoMe99}, \cite{Sh02}). Dans ce
travail nous allons \'{e}tudier des op\'{e}rateurs de Schr\"{o}dinger $%
\mathcal{PT}$-sym\'{e}triques sur l'axe r\'{e}el de la forme:
\begin{equation*}
P=-h^{2}\left( \frac{d}{dx}\right) ^{2}+V(x),
\end{equation*}%
et nous nous placerons dans la limite semi-classique $0<h\rightarrow 0$. La $%
\mathcal{PT}$-symetrie de $\mathcal{P}$ signifie que
\begin{equation}
P\mathcal{PT}=\mathcal{PT}P,  \label{intr.1}
\end{equation}%
Ici les op\'{e}rateurs de parit\'{e} $\mathcal{P}$ et de renversement de
temps $\mathcal{T}$ sont d\'{e}finis par $\mathcal{P}u(x)=u(-x)$, $\mathcal{T%
}u(x)=\overline{u(x)}$, $u$ $\in $ $L^{2}(R)$. Remarquons que (\ref{intr.1})
revient \`{a} la condition suivante sur le potentiel complexe $V$ :
\begin{equation*}
V(-x)=\overline{V(x),}\text{ }x\in \mathbb{R}.
\end{equation*}

Dans le cas d'un op\'{e}rateur de Schr\"{o}dinger auto-adjoint (\`{a}
potentiel r\'{e}el) de la forme
\begin{equation*}
P_{0}=-h^{2}(\frac{d}{dx})^{2}+V_{0}(x),
\end{equation*}%
la $\mathcal{PT}$-sym\'{e}trie equivaut \`a la propriet\'{e} que $V_{0}$ est
pair:
\begin{equation*}
V_{0}(-x)=V_{0}(x).
\end{equation*}

Nous allons \'{e}tudier les cas o\`{u} $V_{0}$ a un puits simple pour un
niveau d'\'{e}nergie donn\'{e} $E_{0}$, et o\`u $V=V_0+i\varepsilon W$ est
une petite perturbation $\mathcal{PT}$-sym\'etrique de $V_0$. Pour cela nous
allons utiliser la m\'{e}thode BKW-complexe (en supposant que $V_{0}$ et $W$
sont analytiques dans un domaine convenable) et \'etablir une condition de
Bohr-Sommerfeld.

Passons maintenant \`{a} la formulation plus pr\'{e}cise de notre r\'{e}%
sultat. Soit $E_{0}\in \mathbb{R}$ un niveau d'\'{e}nergie fix\'{e}. Soit
\begin{equation*}
V_{0}\in C^{\infty }(\mathbb{R}).
\end{equation*}

On fait les hypoth\`eses suivantes sur $V_0$:

\begin{itemize}
\item[\textbf{(H1)}] Il existe un $m_0\ge 0$ tel que pour tout $\alpha \in
\mathbb{N}$, il existe $C_{\alpha }>0$ tel que $|\partial _{x}^{\alpha
}V_{0}(x)|\leq C_{\alpha }(1+|x|)^{m_{0}-\alpha },\ \forall x\in \mathbb{R}$.

\begin{itemize}
\item[(H1.1)] Dans le cas o\`{u} $m_{0}>0$, il existe $C_{0}>0$, tel que $%
V_{0}(x)\geq \frac{1}{C_{0}}|x|^{m_{0}}$, pour $|x|\geq C_{0}$,

\item[(H1.2)] Dans le cas o\`{u} $m_{0}=0,$ $\exists C_{0}>0$, tel que $%
V_{0}(x)\geq E_{0}+\frac{1}{C_{0}}$, quand $|x|\geq C_{0}$.
\end{itemize}

\item[\textbf{(H2)}] $V_{0}$ poss\'{e}de exactement un seul puits de
potentiel:
\begin{equation*}
\begin{split}
\{x\in \mathbb{R};\,V_{0}(x)\leq E_{0}\} =[\alpha _{0}^{0},\beta _{0}^{0}],
\ \{x\in \mathbb{R};\,V_{0}(x)<E_{0}\} =]\alpha _{0}^{0},\beta _{0}^{0}[,
\end{split}%
\end{equation*}%
o\`{u} $-\infty <\alpha _{0}^{0}<\beta _{0}^{0}<+\infty $. De plus, $%
V_{0}^{\prime }(\alpha _{0}^{0})<0,V_{0}^{\prime }(\beta _{0}^{0})>0.$

\item[\textbf{(H3)}] $V_{0}$ admet un prolongement holomorphe \`{a} un
voisinage $U$ dans $\mathbb{C}$ de $\{x\in \mathbb{R};\,V_{0}(x)\leq
E_{0}\}. $
\end{itemize}

Nous allons \'etudier des petites perturbations $\mathcal{PT} $-sym\'{e}%
triques de $P_0$ de la forme
\begin{equation*}
P_{\varepsilon }=-h^{2}(\frac{d}{dx})^{2}+V_{\varepsilon }(x),\
0<\varepsilon \leq 1,
\end{equation*}
o\`{u}
\begin{equation*}
V_{\varepsilon }(x)=V_{0}(x)+i\varepsilon W(x).
\end{equation*}
Ici $W$ est r\'{e}el et impair,
\begin{equation*}
W(-x)=-W(x).
\end{equation*}

Ainsi, on a bien $V_{\varepsilon }(-x)=\overline{V}_{\varepsilon }(x)$ et $%
P_{\varepsilon }$ est $\mathcal{PT}$-sym\'{e}trique. On suppose

\begin{itemize}
\item[\textbf{(H4)}] $W\in C^{\infty }(\mathbb{R},\mathbb{R})$.

\item[\textbf{(H5)}] Comme dans (H1) et avec le m\^{e}me $m_{0}$, pour tout $%
\alpha \in \mathbb{N}$, il existe $C_{\alpha }>0$ tel que $|\partial
_{x}^{\alpha }V_{0}(x)|\leq C_{\alpha }(1+|x|)^{m_{0}-\alpha },\ \forall
x\in \mathbb{R}$.
\begin{equation*}
\forall \alpha \in \mathbb{N},\ \exists C_{\alpha }\text{ telle que }\
|\partial _{x}^{\alpha }W(x)|\leq C_{\alpha }(1+|x|)^{m_{0}-\alpha },\forall
x\in \mathbb{R}.
\end{equation*}

\item[\textbf{(H6)}] $W$ admet un prolongement holomorphe \`{a} un voisinage
$U$ dans $\mathbb{C}$ de $\{x\in \mathbb{R};\,V_{0}(x)\leq E_{0}\}.$
\end{itemize}

\begin{dref}
\label{intr1} Si $U\subset \mathbb{C}$ est un ouvert (comme par exemple dans
les hypoth\`eses ci-dessus), $V$ une fonction holomorphe sur $U$ et $%
P=-h^2(d/dx)^2+V(z)$, on dira que $\alpha =\alpha (E)\in U$ est un point
tournant pour l'\'{e}quation $Pu=Eu$ si $V(\alpha )=E$. Si de plus $%
V^{\prime}(\alpha )\neq 0,$ on dira que $\alpha $\ est un point tournant
simple.
\end{dref}

Soit $D(E_0,\varepsilon )$ le disque ouvert dans $\mathbb{C}$ de centre $E_0$
et de rayon $\varepsilon $. Par le th\'eor\`eme des fonctions implicites on
a:

\begin{prop}
\label{f.m.h} On suppose (H1)--(H6). Il existe $\varepsilon _{0}>0$ telle
que pour $E\in D(E_{0},\varepsilon _{0})$ et $\varepsilon \in
D(0,\varepsilon _{0})$, l'\'{e}quation $V_{\varepsilon }(x)=E$ poss\'{e}de
deux solutions $\alpha _{0}(E,\varepsilon ),$ et $\beta _{0}(E,\varepsilon )$
d\'{e}pendent holomorphiquement de $E$ et $\varepsilon $ avec $\alpha
_{0}(E_{0},0)=\alpha _{0}^{0}$, $\beta _{0}(E_{0},0)=\beta _{0}^{0}$. Ce
sont des points tournants simples pour $P_{\varepsilon
}=(-hd/dx)^{2}+V_{\varepsilon }$.
\end{prop}

Pour $(E,\varepsilon )\in D(E_{0},\varepsilon _{0})\times D(0,\varepsilon
_{0})$ le segment $]\alpha _{0}(E,\varepsilon ),\beta _{0}(E,\varepsilon )[$
appartient a $U$ et $E-V_{\varepsilon }(x)$ ne s'y annule pas. On peut alors
d\'{e}finir la branche continue de la racine carr\'{e}e, $(E$ $%
-V_{\varepsilon }(x))^{\frac{1}{2}}$ pour $(E,\varepsilon )\in
D(E_{0},\varepsilon _{0})\times D(0,\varepsilon _{0})$, $z\in ]\alpha
_{0}(E,\varepsilon ),\beta _{0}(E,\varepsilon )[$ qui est $>0$ quand $E$ est
r\'{e}el et $\varepsilon =0 $. Introduisons l'action
\begin{equation*}
I(E,\varepsilon )=2\int\limits_{\alpha _{0}(E,\varepsilon )}^{\beta
_{0}(E,\varepsilon )}(E-V_{\varepsilon }(z))^{\frac{1}{2}}dz.
\end{equation*}
Ici on int\'{e}gre le long du segment orient\'e qui relie $\alpha
_{0}(E,\varepsilon )$ \`{a} $\beta _{0}(E,\varepsilon )$

\begin{prop}
\label{intr3} Sous les hypoth\`eses (H1)--(H6), si $\ \varepsilon _{0}>0$
est assez petit alors l'action $I(E,\varepsilon )$ est une fonction
holomorphe de $(E,\varepsilon )\in D(E_{0},\varepsilon _{0})\times
D(0,\varepsilon _{0}),$ telle que
\begin{equation*}
I(\overline{E},\varepsilon )=\overline{I(E,\varepsilon )},\text{ quand \ }%
\varepsilon \geq 0.
\end{equation*}
De plus $\dfrac{\partial }{\partial E}I(E,\varepsilon )\neq 0.$
\end{prop}

La preuve sera par des calculs directes en utilisant que $\overline{%
V_{\varepsilon }(-\overline{x})}=V_{\varepsilon }(x)$. Il y aura une preuve
indirecte plus loin.

Soit
\begin{equation*}
P_{\varepsilon }=h^{2}D_{x}^{2}+V_{\varepsilon }(x),\text{ }D_{x}=\dfrac{1}{i%
}\dfrac{d}{dx}
\end{equation*}%
l'op\'{e}rateur de Schr\"{o}dinger sur $\mathbb{R}$, r\'{e}alis\'{e} comme
un op\'{e}rateur ferm\'{e} non born\'{e} $L^{2}(\mathbb{R})\longrightarrow
L^{2}(\mathbb{R})$ de domaine
\begin{equation*}
\mathcal{D}(P_{\varepsilon })=\{u\in L^{2}(\mathbb{R});\,u^{\prime
},u^{\prime \prime },\langle x\rangle ^{m_{0}}u\in L^{2}(\mathbb{R})\},\
\langle x\rangle =(1+x^{2})^{\frac{1}{2}}
\end{equation*}%
On sait alors que le spectre de $P_{\varepsilon }$ dans $D(E_{0},\varepsilon
_{0})$ est discret pour $\varepsilon \in D(0,\varepsilon _{0})$, si $%
\varepsilon _{0}>0$ est assez petit.

Dans le cas $\varepsilon =0$, $P_{0}$ est auto-adjoint, donc les valeurs
propres dans $D(E_{0},\varepsilon _{0})$ sont r\'{e}elles et on a m\^eme,
\begin{equation*}
\inf \sigma _{ess}(P_{0})>E_{0}+\frac{1}{C}.
\end{equation*}
Il est bien connu dans ce cas, que les valeurs propres sont donn\'{e}es par
une condition de quantification de Bohr-Sommerfeld (voir p.ex. \cite{Fe87},
ch II, section 10, \cite{GrSj94}, exercise 12.3).

\begin{theo}
\label{ThBSSp} On fait les hypoth\`{e}ses (H1, H2). Il existe $\varepsilon
_{0},h_{0}>0$ et une fonction r\'{e}elle $\widetilde{I}(E;h)$ de classe $%
C^{\infty }$ sur $]E_{0}-\varepsilon _{0},E_{0}+\varepsilon _{0}[\times
]0,h_{0}[$, admettant un d\'{e}veloppement asymptotique
\begin{equation*}
\widetilde{I}(E;h)\sim I(E)+hI_{1}(E)+...,\ h\rightarrow 0
\end{equation*}%
dans l'espace $C^{\infty }(]E_{0}-\varepsilon _{0},E_{0}+\varepsilon _{0}[)$%
, telle que les valeurs propres de $P_{0}$ dans $]E_{0},-\varepsilon
_{0},E_{0}+\varepsilon _{0}[$ sont donn\'{e}es par la condition de
Bohr--Sommerfeld:
\begin{equation*}
\exists k\in \mathbb{Z},\ E=E_{k},\ \widetilde{I}(E_{k},h)=2k\pi h.
\end{equation*}
\end{theo}

Nous pouvons maintenant \'enoncer le r\'esultat principal de ce travail.

\begin{theo}
\label{ThSP} On fait les hypoth\`{e}ses (H1)--(H6). Il existe $\varepsilon
_{0}>0$ et $h_{0}>0$ tels que $\sigma (P_{\varepsilon })\cap
D(E_{0},\varepsilon _{0})\subset \mathbb{R}$ quand $\ 0\leq \varepsilon \leq
\varepsilon _{0},$ $0<h\leq h_{0}$. Plus pr\'{e}cisement, il existe une
fonction $\widetilde{I}(E,\varepsilon ,h)$ sur $D(E_{0},\varepsilon
_{0})\times D(0,\varepsilon _{0})\times ]0,h_{0}[$, holomorphe en $%
(E,\varepsilon )$, admettant un d\'{e}veloppement asymptotique
\begin{equation*}
\widetilde{I}(E,\varepsilon ;h)\sim I(E,\varepsilon )+hI_{1}(E,\varepsilon
)+...,\ h\rightarrow 0
\end{equation*}%
dans l'espace\ des fonctions holomorphes sur $D(E_{0},\varepsilon
_{0})\times D(0,\varepsilon _{0})$, telle que $\widetilde{I}(E,\varepsilon
;h)\in \mathbb{R}$ quand $E,\varepsilon \in \mathbb{R}$, et telle que pour $%
\varepsilon \in ]0,\varepsilon _{0}[$ les valeurs propres de $P_{\varepsilon
}$ dans $D(E_{0},\varepsilon _{0})$ sont donn\'{e}es par la condition de
Bohr--Sommerfeld:
\begin{equation*}
\exists k\in
\mathbb{Z}
,\ E=E_{k},\ \widetilde{I}(E_{k},\varepsilon ;h)=2k\pi h.
\end{equation*}
\end{theo}

\paragraph{Remerciement.}

Nous tenons \`a remercier Johannes Sj\"ostrand qui nous a propos\'e le sujet
de cette \'etude et qui nous a ensuite soutenu pendant le travail.

\section{M\'ethode BKW complexe en g\'en\'eral}

\label{sw} \setcounter{equation}{0} Dans cette section nous allons revoir
quelques \'{e}l\'{e}ments de la m\'{e}thode BKW complexe, (voir \cite{Ra05},
\cite{Vo81}, \cite{Sj13+}) pour plus de d\'{e}tails. Soit $U\subset \mathbb{C%
}$ un ouvert simplement connexe, d\'{e}signons par $\mathrm{Hol}(U)$
l'espace de Fr\'{e}ch\'{e}t des fonctions holomorphes sur $U$ muni de la
topologie de convergence localement uniforme. Soit $V\in \mathrm{Hol}(U)$ un
potentiel tel que
\begin{equation}
V(z)\neq 0,\ z\in U.  \label{bkw.1}
\end{equation}%
On consid\'{e}re l'\'{e}quation de Schr\"{o}dinger
\begin{equation}
Pu=\left( -h^{2}\left( \frac{d}{dx}\right) ^{2}+V(z)\right) u(z)=0
\label{bkw.2}
\end{equation}%
dans $U$, et on va d'abord chercher une solution BKW formelle de la forme:
\begin{equation}
u(z;h)=a(z;h)e^{i\varphi (z)/h}  \label{bkw.3}
\end{equation}%
o\`{u} $a(z;h)$ a un d\'{e}veloppement asymptotique formel
\begin{equation}
a(z;h)\sim \sum\limits_{j=0}^{+\infty }a_{j}(z)h^{j}  \label{bkw.4}
\end{equation}%
dans l'\'{e}space $\mathrm{Hol}(U)$.

En consid\'{e}rant le d\'{e}veloppement en puissances de $h$ de
\begin{equation*}
e^{-i\varphi (z)/h}(-h^{2}\partial _{z}^{2}+V(z))e^{i\varphi (z)/h }a(z;h)=0,
\end{equation*}%
on trouve
\begin{equation*}
\left( (-(h\partial _z+i\varphi ^{\prime}(z))^2+V(z)\right) a(z;h)=0
\end{equation*}
o\`u plus explicitement,
\begin{equation*}
\left( -(h\partial _{z})^{2}+{\varphi ^{\prime }}(z)^{2}-2i\varphi ^{\prime
}h\partial _{z}-ih\varphi ^{\prime\prime}+V(z)\right) a(z) =0.
\end{equation*}
On est amen\'e \`a choisir $\varphi $ solution de l'\'{e}quation eiconale
\begin{equation}
(\varphi ^{\prime }(z))^{2}+V(z)=0,  \label{bkw.5}
\end{equation}
qu'on peut facilement r\'esoudre sur $U$:

\begin{prop}
\label{bkw1} On suppose qui'il n'y a pas de points tournants c'est \`{a}
dire que $V(z)\neq 0$ pour tous $z\in U$. L'\'{e}quation eiconale poss\'{e}%
de deux solutions holomorphes \`{a} des constantes pr\'{e}s, donn\'{e}es
par:
\begin{equation}  \label{bkw.6}
\varphi (z)=\pm \int_{z_{0}}^{z}(-V(w))^{\frac{1}{2}}dw,
\end{equation}
o\`u $z_0\in U$. Ici, $(-V(w))^{\frac{1}{2}}$ d\'{e}signe une branche
holomorphe de la racine carr\'e de $-V(z)$ sur $U$.
\end{prop}

Il reste ensuite \`a chercher un d\'eveloppement formel, (\ref{bkw.4}), tel
que
\begin{equation}  \label{bkw.7}
\left( \partial _{z}\varphi (z)\partial _{z}+\frac{\partial _{z}^{2}\varphi
(z)}{2} -i\frac{h\partial _{z}^{2}}{2}\right) a(z;h) =0.
\end{equation}
En annulant successivement les puissances de $h$, on trouve une suite d'\'{e}%
quations de transport:
\begin{equation}  \label{bkw.8}
\begin{split}
(\partial _{z}\varphi \, \partial _{z}+\frac{\partial _{z}^{2}\varphi }{2}
)a_{0} &=0, \\
(\partial _{z}\varphi \, \partial _{z}+\frac{\partial _{z}^{2}\varphi }{2}
)a_{1} &=i\frac{\partial _{z}^{2}}{2}a_{0}, \\
(\partial _{z}\varphi \, \partial _{z}+\frac{\partial _{z}^{2}\varphi }{2})
a_{k} &=i\frac{\partial _{z}^{2}}{2}a_{k-1},\ k\geq 1 .
\end{split}%
\end{equation}
La solution de la premi\`ere \'equation de transport est donn\'ee par
\begin{equation*}
a_0(z)=C(\partial _z\varphi )^{-\frac{1}{2}}=C(-V(z))^{-\frac{1}{4}}.
\end{equation*}

\begin{prop}
\label{bkw2} Soit $z_0\in U$ et fixons une solution de l'\'equation eiconale
(\ref{bkw.5}). Soient $a_{0}^{0},\, a_{1}^{0},\, a_{2}^{0},\,...$ des
nombres complexes arbitraires. Alors il existe une unique solution BKW
formelle de l'\'equation (\ref{bkw.2}) de la forme (\ref{bkw.3}), (\ref%
{bkw.4}) avec
\begin{equation*}
a_{0}(z_0)=a_{0}^{0},\ a_{1}(z_�%
{{}^\circ}%
)=a_{1}^{0},\ a_{2}(z_0)=a_{2}^{0},...
\end{equation*}
\end{prop}

\begin{dref}
\label{bkw3}On appelle ligne de Stokes une courbe $\gamma :[ a,b]
\rightarrow U$ de classe $C^{1}$, telle que
\begin{equation*}
{\Im}m\int_{s}^{t}(-V(\gamma (\tau ))^{\frac{1}{2}}d\gamma (\tau )=0
\end{equation*}
pour tous $s,t$ $\in [ a,b ] $. Autrement dit, ${\Im}m \varphi $ doit
\^etre constant sur toute ligne de Stokes.
\end{dref}

Le r\'esultat suivant permet de passer des solutions formelles \`a des
solutions exactes en respectant la r\`egle fondamentale de la m\'ethode BKW
complexe qui est de se d\'eplacer toujours dans la direction o\`u le facteur
phase $\exp (i\varphi /h)$ est croissant en module, donc en particulier
transversalement aux lignes de Stokes.

\begin{theo}
\label{bkw4} Soit $-\infty <b<c<+\infty $ et $\gamma :[b,c]\rightarrow U$
une courbe de classe $C^{1}$, telle que $\frac{d}{dt}(-{\Im}m\varphi
(\gamma (t))>0$, $b\leq t\leq c$. Soit
\begin{equation}
u_{\mathrm{BKW}}\sim (a_{0}(z)+ha_{1}(z)+...)e^{i\varphi (z)/h}
\label{bkw.9}
\end{equation}%
une solution BKW formelle de (\ref{bkw.2}).

\begin{enumerate}
\item Il existe une solution exacte $u$ de (\ref{bkw.2}) et un voisinage
ouvert $B$ de $\gamma (b)$ tels que
\begin{equation}  \label{bkw.10}
\begin{split}
u(z;h)&=a(z;h)e^{i\varphi (z)/h}\text{ dans }B, \\
a(z;h)&\sim a_0(z)+ha_1(z)+... \text{ dans }\mathrm{Hol\,}(B).
\end{split}%
\end{equation}
$B$ ne d\'epend pas du choix de la solution BKW formelle.

\item Il existe un voisinage ouvert $\Gamma $

\item de $\gamma (]b,c])$ tel que si $u$ est une solution exacte comme dans
1), alors la description (\ref{bkw.10}) s'\'{e}tend \`{a} $\Gamma $.
\end{enumerate}
\end{theo}

Ce r\'esultat est bien connu. Voir par exemple \cite{Sj13+}.

On retourne maintenant \`a la situation d\'ecrite dans la section \ref{intr}
et on adopte les hypoth\`eses (H1)--(H6). Nous avons d\'ej\`a d\'efini le
puits $[\alpha _0^0,\beta _0^0]$ pour $V_0-E_0$ et les points tournants $%
\alpha _0(E,{\varepsilon})$, $\beta _0(E,{\varepsilon})$.

On s'int\'eresse aux solutions \`a d\'ecroissance exponentielle pr\`es de $%
\pm \infty $ de l'\'equation
\begin{equation}  \label{bkw.11}
(-(h\partial _x)^2+V_{{\varepsilon}} (x)-E)u(x)=0.
\end{equation}
Introduisons les espaces vectoriels complexes
\begin{equation*}
\mathcal{E}_{\pm}=\mathcal{E}_{\pm}(E,{\varepsilon};h)=\{ u\in C^\infty (%
\mathbb{R});\ u\text{ v\'erifie (\ref{bkw.11}) et }u \text{ est born\'e sur }%
\mathbb{R}_{\pm}\}.
\end{equation*}
Gr\^ace \`a l'ellipticit\'e de $P_{\varepsilon}-E$ pr\`es de $\pm \infty $,
on a le r\'esultat bien connu suivant:

\begin{prop}
\label{bkw5} $\mathrm{dim\,}\mathcal{E}_{\pm}=1$, c'est \`a dire chaque
espace est engendr\'e par une seule solution de (\ref{bkw.11}): $\mathcal{\ E%
}_{\pm}=\mathbb{C}u_{\pm}$. La fonction $u_{\pm}$ est \`a d\'ecroissance
exponentielle pr\`es de $x=\pm \infty $.
\end{prop}

On peut aussi d\'ecrire le comportement asymptotique de $u_{\pm}$ pr\`es de $%
\pm \infty $. Commen\c cons par employer la m\'ethode BKW formelle sur $%
]-\infty ,\alpha _0^0-\delta _0]$ et sur $[\beta _0^0+\delta _0,+\infty [$
quand $\delta _0>0$ et pour ${\varepsilon}_0>0$ assez petit en fonction $%
\delta _0$. L'analyse sur les deux intervalles est essentiellement la m\^eme
et on va se concentrer sur $[\beta _0^0+\delta _0,+\infty [$. L'\'equation
eiconale
\begin{equation}  \label{bkw.12}
(\varphi ^{\prime}(x))^2+V_{\varepsilon}(x)-E=0
\end{equation}
poss\`ede la solution
\begin{equation}  \label{bkw.13}
\varphi (x)=i\int_{\beta _0(E,\varepsilon )}^x (V_{{\varepsilon}%
}(y)-E)^{1/2} dy,\ x\ge \beta _0^0+\delta _0,
\end{equation}
o\`u on choisit la branche de la racine carr\'ee qui d\'epend continuement
de $(E,{\varepsilon})$ et qui est $>0$ quand $(E,{\varepsilon})=(E_0,0)$. Il
est alors clair que
\begin{equation}  \label{bkw.14}
\partial _x^{\alpha} \varphi (x)=\mathcal{O}(1)(1+|x|)^{\frac{m_0}{2}%
+1-\alpha },\ \alpha \in \mathbb{N}.
\end{equation}

On cherche ensuite une solution BKW formelle comme dans (\ref{bkw.3}), (\ref%
{bkw.4}). Alors le symbole $a$ doit v\'erifier (\ref{bkw.7}) c.\`a.d. la
suite des \'equations de transport (\ref{bkw.8}). Ici on peut prendre $%
a_0(x)=(\partial _x\varphi )^{-1/2}$ et si on pose $a_k(x)=f_k(x)a_0(x)$ on
trouve $f_0=1$ et
\begin{equation}  \label{bkw.15}
\partial _xf_k=\frac{i}{2\partial _x\varphi }\left(\partial _x^2f_{k-1}+2%
\frac{\partial _x a_0}{a_0}\partial _xf_{k-1}+\frac{\partial _x^2a_0}{a_0}%
f_{k-1} \right),\ k\ge 1.
\end{equation}
Par r\'ecurrence sur $k$ on voit qu'on peut trouver des solutions $f_1$, $%
f_2 $, ..., tels que
\begin{equation}  \label{bkw.16}
\partial _x^{\alpha} f_k=\mathcal{O}(1)(1+|x|)^{-k(1+m_0/2)-\alpha },
\end{equation}
donc pour les $a_k$,
\begin{equation}  \label{bkw.17}
\partial _x^{\alpha} a_k(x)=\mathcal{O}(1)(1+|x|)^{-m_0/4-k(1+m_0/2)-\alpha
}.
\end{equation}

Par des arguments standard d'\'equations diff\'erentielles ordinaires on
peut ensuite passer des solutions formelles aux solutions exactes pour
arriver \`a:

\begin{prop}
\label{bkw6} $\forall \delta _{0}>0$, $\exists \varepsilon_0 >0$ tel que
pour $(E,\varepsilon )\in ((E_{0},\varepsilon _{0})\times D(0,\varepsilon
_{0}),$ l'\'{e}quation (\ref{bkw.11}) ait une solution holomorphe en $%
(E,\varepsilon )$ de la forme
\begin{equation*}
u_+(z;h)=a(z;h)e^{i\varphi (z)/h}\text{ o\`{u} }a(z;h)\sim
\sum\limits_{j=0}^{\infty }a_{j}h^{j}\text{ sur }[\beta _0^0+\delta
_0,+\infty [
\end{equation*}
au sens suivant: Pour tout $(N,\alpha )\in \mathbb{N}^*\times \mathbb{N}$ il
existe une constante $C_{N,\alpha }>0$ telle que
\begin{equation*}
| \partial ^{\alpha }(a(x;h)-\sum\limits_{j=0}^{N-1}a_{k}(x)h^{k}| \leq
C_{N,\alpha }h^{N}(1+\mid x\mid )^{-\frac{m_{0}}{4}-N(\frac{m_{0}}{2}
+1)-\alpha }
\end{equation*}
pour $x\in [\beta _0^0+\delta _0,+\infty [$. Ici $a_0=(\partial _x\varphi
)^{-1/2}$ et $a_k$ v\'erifie (\ref{bkw.17}).
\end{prop}

\begin{remark}
Nous avons le m\^{e}me r\'{e}sultat dans un intervalle $] -\infty ,\alpha
_{0}^{0},-\delta _{0}] $ o\`{u} on choisit la branche oppos\'{e} de la
racine carr\'{e}e $(V_{{\varepsilon}}(x)-E)^{\frac{1}{2}}$.
\end{remark}

\section{Analyse BKW pr\`es d'un point tournant simple}

\label{pts} \setcounter{equation}{0} Dans cette section on suit la pr\'{e}%
sentation dans \cite{Sj13+} de pr\`{e}s. Soit $\Omega \subset \mathbb{C}$ un
ouvert simplement connexe, $V\in \mathrm{Hol\,}(\Omega )$. Soit $z_{0}\in
\Omega $ un point tournant simple, {%
\begin{equation}
V(z_{0})=0,\ V^{\prime }(z_{0})\neq 0.  \label{pts.1}
\end{equation}%
} Pour simplifier la notation on suppose que $z_{0}=0$ et on s'int\'{e}resse
aux solutions de (\ref{bkw.2}) qui dans certaines r\'{e}gions prennent la
forme $a(z;h)e^{\varphi (z)/h}$ (sans facteur $i$ dans l'exposant pour
simplifier les notations). Consid\'{e}rons l'\'{e}quation eiconale, {%
\begin{equation}
\varphi ^{\prime }(z)=V(z)^{\frac{1}{2}}  \label{pts.2}
\end{equation}%
} dans un voisinage de $0$. (On diminuera $\Omega $ autour de $z=0$ chaque
fois que cela nous arrange). Il est clair que $\varphi (z)$ sera multi-valu%
\'{e} en g\'{e}n\'{e}ral et pour mieux comprendre la structure de cette
singularit\'{e} on passe au recouvrement double de $\Omega \setminus \{0\}$,
en posant $z=w^{2}$. Alors
\begin{equation*}
\frac{\partial }{\partial z}=\frac{1}{2w}\frac{\partial }{\partial w},
\end{equation*}%
et si on pose $\widetilde{V}(w)=V(z)=F(z)z=F(w^{2})w^{2}$, $\varphi (z)=%
\widetilde{\varphi }(w)$, o\`{u} $F(0)\neq 0$, l'\'{e}quation eiconale
devient
\begin{equation*}
\partial _{w}\widetilde{\varphi }=F(w^{2})^{\frac{1}{2}}2w^{2},
\end{equation*}%
o\`{u} le membre droit est une fonction holomorphe paire. Si on exige aussi
que $\varphi (0)=\widetilde{\varphi }(0)=0$, on voit que $\widetilde{\varphi
}(w)$ est une fonction holomorphe impaire de la forme
\begin{equation*}
\widetilde{\varphi }(w)=\frac{2}{3}\widetilde{F}(w^{2})w^{3},\hbox{ o\`u }%
\widetilde{F}(0)=F(0)^{\frac{1}{2}}=V^{\prime }(0)^{\frac{1}{2}}.
\end{equation*}%
Dans le coordonn\'{e} $z$ on obtient une fonction double-valu\'{e}e, {%
\begin{equation}
\varphi (z)=\frac{2}{3}\widetilde{F}(z)z^{\frac{3}{2}}.  \label{pts.3}
\end{equation}%
}

Cherchons maintenant des lignes de Stokes et anti-Stokes qui passent par $0$%
. (Comme on a supprim\'{e} le facteur $i$ dans l'exposant dans les repr\'{e}%
sentations BKW, ${\Re}e\varphi =\mathrm{Const.}$ sur chaque ligne de
Stokes et par d\'{e}finition ${\Im}m\varphi =\mathrm{Const}$ sur les
lignes anti-Stokes). Sur de telles courbes nous avons ${\Re}e\varphi =0$ o%
\`{u} ${\Im}m\varphi =0$, c'est \`{a} dire ${\Im}m\varphi ^{2}=0$: $%
{\Im}m\widetilde{F}(z)^{2}z^{3}=0$. Autrement dit, $\widetilde{F}%
(z)^{2}z^{3}=t^{3}$ pour un $t\in \mathrm{vois\,}(0,\mathbb{R})$ et en
prenant la racine cubique nous obtenons trois courbes $\gamma _{k}$
\begin{equation*}
\widetilde{F}(z)^{\frac{2}{3}}z=e^{2\pi ik/3}t,\ k\in \{0,1,2\}\simeq
\mathbf{Z}/3\mathbf{Z}.
\end{equation*}%
On obtient la figure suivante o\`{u} on a pris $V^{\prime }(0)>0$ pour fixer
les id\'{e}es:
Chaque courbe $\gamma _{k}\setminus \{0\}$ se d\'{e}compose en une ligne de
Stokes $\gamma _{k}^{-}$ et une ligne d'anti-Stokes $\gamma _{k}^{+}$. Les
trois lignes de Stokes et le point tournant d\'{e}limitent trois
\textquotedblleft secteurs de Stokes\textquotedblright\ ferm\'{e}s $\Sigma
_{k}$. Dans la figure $1$ nous avons aussi trac\'{e} quelques
lignes de Stokes \`{a} l'int\'{e}rieur de chaque secteur.

\begin{figure}[!h]
 \centering
 {\includegraphics[height=8cm,width=10cm]{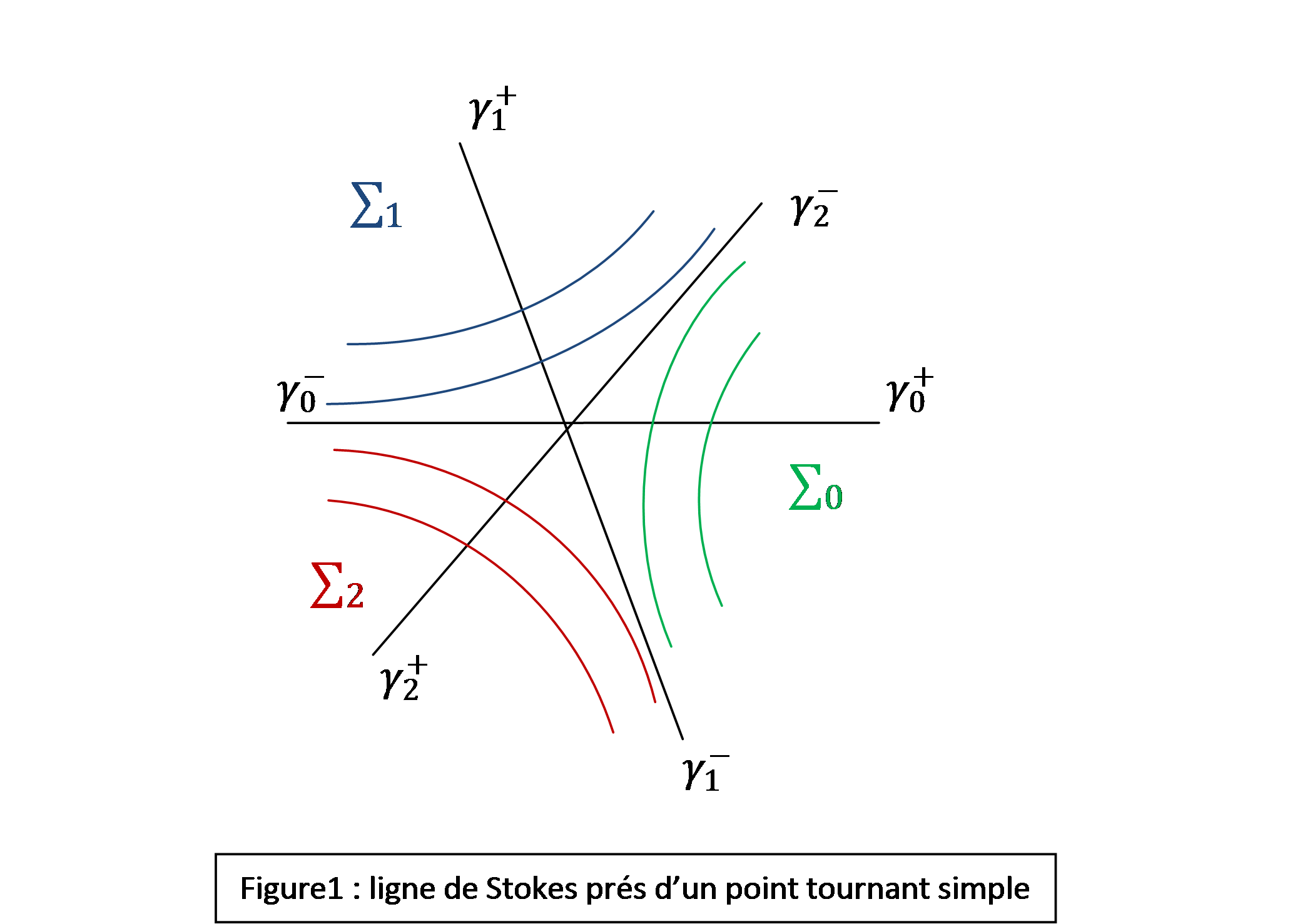}\label{Fig1n}}
\caption{Lignes de Stokes pr\'es d'un point tournant simple.}
\end{figure}

Soit $\varphi _k$ la branche de $\varphi $ dans $\Omega \setminus \gamma
_k^- $ telle que ${\Re}e \varphi _k<0$ dans $\mathrm{int\,}(\Sigma _k)$, $%
\varphi _k(0)=0$. Remarquons que $\varphi _{k+1}$ et $\varphi _k$ sont tous
les deux bien d\'efinis dans $\Sigma _k\cup \Sigma _{k+1}$ et y v\'erifient $%
\varphi _{k+1}=-\varphi _k$.

D'apr\`es le principe fondamental de la m\'ethode BKW complexe il existe des
solutions exactes $u=u_j$, $j\in \mathbf{Z}/3\mathbf{Z}$ de l'\'equation $%
(-(h\partial )^2+V)u=0$ dans $\Omega $ telles que
\begin{equation}  \label{pts.3.2}
\begin{cases}
u_j(z;h)=a_j(z;h)e^{\varphi _j(z)/h} \\
a_j(z;h)\sim a_{j,0}(z)+ha_{j,1}(z)+...%
\end{cases}
\hbox{ dans } \mathrm{int\,}(\Sigma _j).
\end{equation}
Cette description asymptotique s'\'etend au complement d'un voisinage
arbitrairement petit de $\gamma _j^-\cup \{0 \}$ (qui peut \^etre atteint de
$\Sigma _j$ par des chemins qui ne sont jamais tangents aux lignes de
Stokes). On rappelle aussi que $a_{j,0}$ est unique \`a un facteur constant
pr\`es et qu'on peut choisir
\begin{equation}  \label{pts.3.4}
a_{j,0}(z)=(\varphi ^{\prime}_j(z))^{-\frac{1}{2}},
\end{equation}
o\`u pour l'instant on ne fixe pas la branche de la racine carr\'ee.

Rappelons que si $u,v$ sont des solutions de notre \'{e}quation de Schr\"{o}%
dinger homog\`{e}ne, alors le Wronskien
\begin{equation*}
W(u,v)=(h\partial u)v-u(h\partial v).
\end{equation*}%
est constant. Appliquant l'asymptotique de $u_{0}$ et $u_{1}$ en un point de
$\mathrm{int\,}(\Sigma _{0}\cup \Sigma _{1})$, nous voyons que $%
W(u_{0},u_{1})$ a un d\'{e}veloppement asymptotique en puissances de $h$:
\begin{equation*}
\begin{split}
W(u_{0},u_{1})& =2a_{0,0}a_{1,0}\partial \varphi _{0}+\mathcal{O}(h) \\
& =\frac{2\varphi _{0}^{\prime }}{\sqrt{\varphi _{0}^{\prime }\varphi
_{1}^{\prime }}}+\mathcal{O}(h) \\
& =\frac{2\varphi _{0}^{\prime }}{\sqrt{-(\varphi _{0}^{\prime })^{2}}}+%
\mathcal{O}(h) \\
& =\pm 2i+\mathcal{O}(h).
\end{split}%
\end{equation*}%
De la m\^{e}me fa\c{c}on
\begin{equation*}
\begin{split}
W(u_{1},u_{2})& =\pm 2i+\mathcal{O}(h) \\
W(u_{2},u_{0})& =\pm 2i+\mathcal{O}(h).
\end{split}%
\end{equation*}

On peut d\'{e}terminer les signes de la fa\c{c}on suivante: Fixons une
branche de $(\varphi _{j}^{\prime })^{1/2}$ comme ci-dessus pour $j=0,1,2\
\mathrm{mod\,}4\mathbf{Z}$. Alors pour deux secteurs de Stokes diff\'{e}%
rents, $j\neq k$ nous avons dans l'int\'{e}rieur de $\Sigma _{j}\cup \Sigma
_{k}$ que
\begin{equation}
(\varphi _{j}^{\prime })^{1/2}=i^{\nu _{j,k}}(\varphi _{k}^{\prime })^{1/2},
\label{pts.3.5}
\end{equation}%
o\`{u} $\nu _{j,k}\in \mathbf{Z}/4\mathbf{Z}$ est impair et $\nu _{j,k}=-\nu
_{k,j}$.

En commen\c{c}ant dans $\Sigma _{0}$ on fait un tour autour de $0$ dans le
sens positif et on note que
\begin{equation}
\begin{split}
(\varphi _{1}^{\prime })^{1/2}& =i^{\nu _{1,0}}(\varphi _{0}^{\prime })^{1/2}
\\
(\varphi _{2}^{\prime })^{1/2}& =i^{\nu _{2,1}}(\varphi _{1}^{\prime })^{1/2}
\\
(\varphi _{0}^{\prime })^{1/2}& =i^{\nu _{0,2}}(\varphi _{2}^{\prime
})^{1/2}.
\end{split}
\label{pts.3.6}
\end{equation}%
Cela veut dire que si on suit une branche continue de $(\varphi _{0}^{\prime
})^{1/2}$ autour de $0$ dans le sens positif, alors apr\`{e}s un tour, on
obtient la branche
\begin{equation*}
i^{-(\nu _{0,2}+\nu _{2,1}+\nu _{1,0})}(\varphi _{0}^{\prime })^{1/2}.
\end{equation*}%
Mais $\varphi _{0}^{\prime }=V^{1/4}$ pour une branche convenable de la
racine quatri\`{e}me et si on suit cette fonction autour de $0$ une fois, on
trouve $iV^{1/4}$. Ceci donne la condition de co-cycle,
\begin{equation}
\nu _{0,2}+\nu _{2,1}+\nu _{1,0}\equiv -1\ \mathrm{mod\,}4\mathbf{Z}.
\label{pts.4}
\end{equation}%
Nous pouvons maintenant pr\'{e}ciser les signes dans les calculs des
Wronskiens ci-dessus:
\begin{equation}
W(u_{j},u_{k})=\frac{-2\varphi _{k}^{\prime }}{\sqrt{\varphi _{j}^{\prime }}%
\sqrt{\varphi _{k}^{\prime }}}+\mathcal{O}(h)=\frac{-2\varphi _{k}^{\prime }%
}{i^{\nu _{j,k}}\sqrt{\varphi _{k}^{\prime }}^{2}}+\mathcal{O}(h)=2i^{\nu
_{j,k}}+\mathcal{O}(h).  \label{pts.5}
\end{equation}

L'espace de solutions nulles est de dimension 2, donc nous avons une
relation {%
\begin{equation}
\alpha _{-1}u_{-1}+\alpha _{0}u_{0}+\alpha _{1}u_{1}=0,  \label{pts.6}
\end{equation}%
} o\`{u} le vecteur $(\alpha _{-1},\alpha _{0},\alpha _{1})^{\mathrm{t}}\in
\mathbf{C}^{3}\setminus \{0\}$ est bien d\'{e}fini \`{a} un facteur scalaire
pr\`{e}s. Appliquant $W(u_{j},\cdot )$ \`{a} cette relation, on obtient {%
\begin{equation}
(W(u_{j},u_{k}))_{j,k}%
\begin{pmatrix}
\alpha _{-1} \\
\alpha _{0} \\
\alpha _{1}%
\end{pmatrix}%
=0,  \label{pts.7}
\end{equation}%
} o\`{u} plus explicitement, {%
\begin{equation}
\begin{pmatrix}
0 & a & b \\
-a & 0 & c \\
-b & -c & 0%
\end{pmatrix}%
\begin{pmatrix}
\alpha _{-1} \\
\alpha _{0} \\
\alpha _{1}%
\end{pmatrix}%
=0.  \label{pts.8}
\end{equation}%
} Nous pouvons prendre {%
\begin{equation}
\begin{pmatrix}
\alpha _{-1} \\
\alpha _{0} \\
\alpha _{1}%
\end{pmatrix}%
=%
\begin{pmatrix}
c \\
-b \\
a%
\end{pmatrix}%
,  \label{pts.9}
\end{equation}%
} donc \`{a} un facteur commun pr\`{e}s, nous avons {%
\begin{equation}
\alpha _{j}=\pm i+\mathcal{O}(h).  \label{pts.10}
\end{equation}%
}

(\ref{pts.5}) permet de pr\'eciser les valeurs de $a,b,c$ et de $\alpha
_{-1},\alpha _0,\alpha _1$:
\begin{equation*}
\begin{split}
a&=\frac{1}{2}W(u_{-1},u_0)=i^{\nu _{-1,0}}+\mathcal{O}(h) \\
b&=\frac{1}{2}W(u_{-1},u_1)=i^{\nu _{-1,1}}+\mathcal{O}(h) \\
c&=\frac{1}{2}W(u_0,u_1)=i^{\nu _{0,1}}+\mathcal{O}(h),
\end{split}%
\end{equation*}
(apr\`es l'insertion d'un facteur commun 1/2) ce qui donne
\begin{equation}  \label{pts.11}
\begin{pmatrix}
\alpha _{-1} \\
\alpha_0 \\
\alpha _1%
\end{pmatrix}%
=
\begin{pmatrix}
i^{\nu _{0,1}} \\
-i^{\nu _{-1,1}} \\
i^{\nu _{-1,0}}%
\end{pmatrix}
+\mathcal{O}(h)=%
\begin{pmatrix}
i^{\nu _{0,1}} \\
i^{\nu _{1,-1}} \\
i^{\nu _{-1,0}}%
\end{pmatrix}
+\mathcal{O}(h).
\end{equation}

\begin{remark}
\label{pts1} Parfois il est plus naturel de changer les notations, en
\'ecrivant $i\varphi _j$ dans (\ref{pts.3.2}) \`a la place de $\varphi _j$
de telle sorte $u_j(z;h)=a_j(z;h)e^{i\varphi _j(z)/h}$ avec $\mathrm{Im\,}
\varphi _j\ge 0$ dans $\Sigma _j$. (\ref{pts.3.4}) devient alors $%
a_{j,0}(z)=(i\varphi _j^{\prime})^{-1/2}=V(z)^{-1/4}$ et dans (\ref{pts.3.5}%
), (\ref{pts.3.6}) on doit remplacer $\varphi_j^{\prime}$ par $i\varphi
_j^{\prime}$.
\end{remark}

\section{Quantification de Bohr Sommerfeld pour un puits de potentiel sans $%
\mathcal{PT}$ sym\'{e}trie}

\label{qbss}%
\setcounter{equation}{0}%
Soit $V_{0}$ un potentiel analytique \`{a} valeurs r\'{e}elles sur un
voisinage r\'eel de $[A,B]$, o\`u $-\infty <A<B<+\infty $. Soit $E_{0}$ $\in
\mathbb{R}$ et supposons qu'il existe $A<\alpha _0^0<\beta _0^0<B$ tels
que
\begin{equation}
V_{0}-E_{0}\
\begin{cases}
>0\text{ dans }[A,\alpha _0^0[\cup ]\beta _0^0,B], \\
<0\text{ dans }]\alpha _0^0,\beta _0^0[.%
\end{cases}
\label{qbss.1}
\end{equation}
On suppose aussi que $\alpha _0^0$ et $\beta _0^0$ sont deux points
tournants simples pour $V_{0}(x)-E_{0}:$

\begin{equation*}
V_{0}^{\prime }(\alpha _{0}^{0})<0,V_{0}^{\prime }(\beta _{0}^{0})>0,
\end{equation*}%
L'une des lignes de Stokes de $\alpha _{0}^{0}$ atteint $\beta _{0}^{0}.$

\begin{figure}[!h]
 \centering
 {\includegraphics[height=5cm,width=10cm]{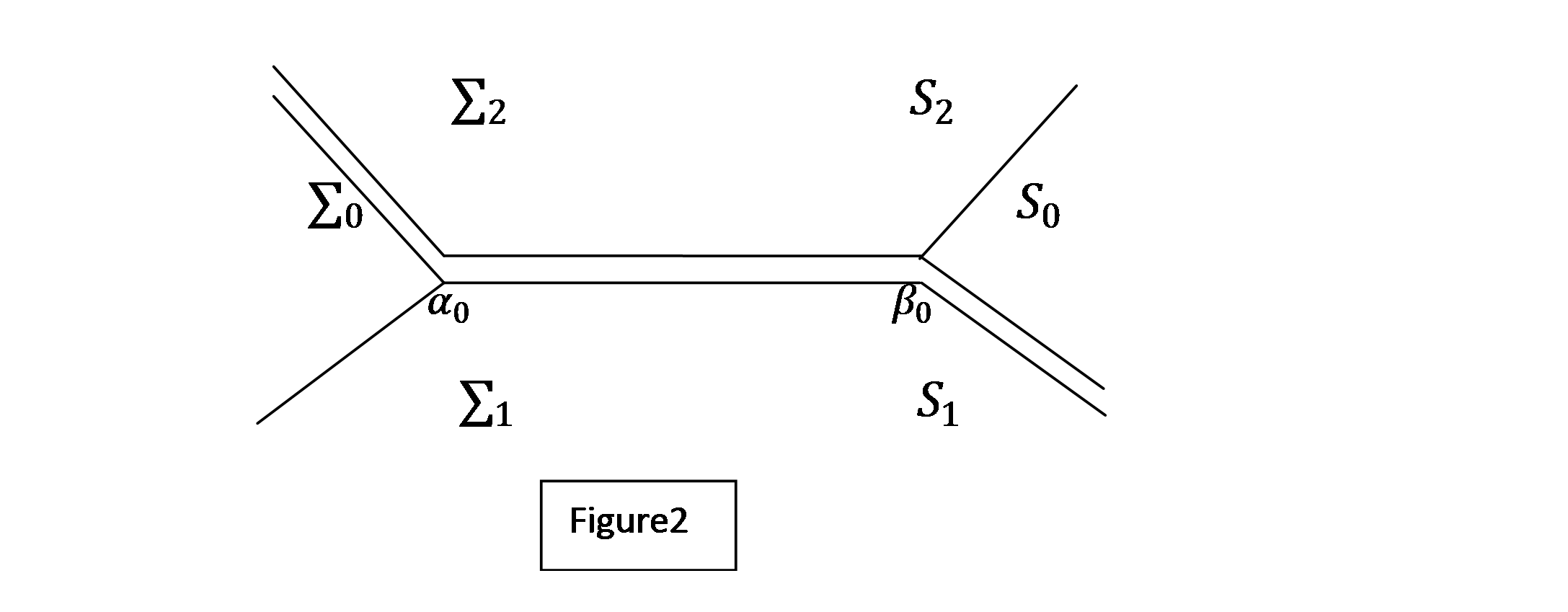}\label{Fig2n}}
\caption{Raccordement pour deux points tournants}
\end{figure}

Soit $U\Subset \mathbb{C} $ un voisinage complexe de $[A,B]$ dans lequel $%
V_0 $ s'\'{e}tend holomorphiquement. Soit $V_{\varepsilon
}(x)=V_{0}(x)+i\varepsilon W(x)$ o\`{u} $W(x)$ est une fonction holomorphe
dans $U$.

Si $\varepsilon \in \mathbb{C}$ est assez petit en module et $E$ appartient
\`{a} un petit voisinage complexe de $E_{0}$, nous avons encore deux points
tournants simples $\alpha _{0}(E,\varepsilon )$, $\beta _{0}(E,\varepsilon )$
proches de $\alpha _{0}^{0}$ et $\beta _{0}^{0}$ qui d\'{e}pendent
holomorphiquement de $(E,\varepsilon )$.


Le dessin indique les trois secteurs de Stokes $\Sigma _{j}$ proches de $%
\alpha _0$ et trois secteurs de Stokes $S_{j}$ proches de $\beta _0$ pour $%
j=-1,0,1$. Notons que $A\in \Sigma _{0}$ et $B\in $ $S_{0}$.

Pour chaque secteur $\Sigma _{j}$ on a une solution exacte $u=u_{j}$ dans $U$
de l'\'equation de Schr\"odinger $(-h^2\partial ^2+V_{\varepsilon} -E)u=0$
telle que
\begin{equation}
u_{j}(z;h)=a_{j,\alpha _0}(z;h)e^{i\varphi _{_{^{j,\alpha _0}}}(z)/h}\text{
dans }\Sigma _{j},  \label{qbss.2}
\end{equation}
avec $\varphi _{j,\alpha _0}(\alpha _0)=0$, $\mathrm{Im\,} \varphi
_{j,\alpha _0}>0$ dans l'int\'erieur de $\Sigma _j$. De m\^{e}me, on a une
solution exacte $v_{j}$ telle que
\begin{equation}
v_{j}(z;h)=a_{j,\beta _0}(z;h)e^{i\varphi _{_{^{j,\beta _0}}}(z)/h}\text{
dans }S _{j},  \label{qbss.3}
\end{equation}
avec $\varphi _{j,\beta _0}(\beta _0)=0$, $\mathrm{Im\,} \varphi _{j,\beta
_0}>0$ dans l'int\'erieur de $S_j$.

Quitte \`a diminuer $\Sigma _0$ et $S_0$ pour que $A\not\in \Sigma _0$, $%
B\not\in S_0$ nous pouvons nous arranger pour que $u_{0}(A)$ $=0$, $%
v_{0}(B)=0$ et que $u_0$, $v_0$ d\'ependent de mani\`ere holomorphe de $%
(E,\varepsilon )$. De la m\^eme fa\c con, pour $j=\pm 1$, on peut s'arranger
pour que $u_j$ et $v_j$ d\'ependent holomorphiquement de $(E,\varepsilon )$.

Maintenant, consid\'{e}rons le probl\`{e}me de Dirichlet
\begin{equation}
-((h\partial )^{2}+V_{\varepsilon }-E)u=0,u(A)=u(B)=0.  \label{qbss.4}
\end{equation}%
En d'autres termes, nous cherchons le spectre de l'op\'{e}rateur non born%
\'{e}
\begin{equation}
P_{\varepsilon }=-(h\partial ^{2})+V_{\varepsilon }:L^{2}(]A,B[)\rightarrow
L^{2}(]A,B[),  \label{qbss.5}
\end{equation}%
de domaine\footnote{%
Le cas d'un op\'{e}rateur d\'{e}fini sur tout l'axe r\'{e}el comme dans les
th\'{e}or\`{e}mes \ref{ThBSSp}, \ref{ThSP} se traite de la m\^{e}me fa\c{c}%
on avec des modifications mineures.}
\begin{equation}
\mathcal{D(}P_{\varepsilon })=\{u\in H^{2}(]A,B[);\,u(A)=u(B)=0\}.
\label{qbss.6}
\end{equation}%
Nous voyons que
\begin{equation}
E\in \sigma (P_{\varepsilon })\Leftrightarrow \mathcal{W(}u_{0},v_{0})=0.
\label{qbss.7}
\end{equation}

Pour $j=\pm 1$ nous pouvons choisir $u_j, v_j$ collin\'eaires:
\begin{equation}  \label{qbss.8}
\begin{split}
u_{j}(z;h)&=c(h)e^{i\varphi _{j,\alpha _0}(\beta _0)/h}v_{j}, \\
c(h)&\sim c_0(E,\varepsilon )+hc_1(E,\varepsilon )+...\text{ dans }\mathrm{%
Hol\,}(\mathrm{vois\,}((E_0,0),\mathbb{C}^2)).
\end{split}%
\end{equation}%
Ici,
\begin{equation}  \label{qbss.9}
\varphi _{j,\alpha _0}(\beta _0)=\pm\int_{\alpha _0}^{\beta
_0}(E-V_{\varepsilon }(z))^{\frac{1}{2}}dz,
\end{equation}
avec le signe $+$ pour $j=-1$ et le signe $-$ pour $j=1$.

Quand $E=E_0$, $\varepsilon =0$, alors pour $z<\alpha _0^0$ (r\'eel),
\begin{equation}  \label{qbss.12}
u_0=\frac{1+\mathcal{O}(h)}{\sqrt{i\varphi _0^{\prime}}}e^{i\varphi
_{0,\alpha _0}/h}
\end{equation}
est \`a valeurs r\'eelles, $i\varphi _{0,\alpha _0}<0$, $i\varphi _{0,\alpha
_0}^{\prime}>0$. Dans la discussion de la section \ref{qbss} on peut choisir
\begin{equation}  \label{qbss.13}
\nu _{1,0}=-1,\ \nu _{0,-1}=-1,\ \nu _{-1,1}=1,
\end{equation}
respectant la condition (\ref{pts.4}). La relation (\ref{pts.11}) devient
\begin{equation}  \label{qbss.14}
\alpha _{-1}=i+\mathcal{O}(h),\ \alpha _0=-i+\mathcal{O}(h),\ \alpha _1=i+%
\mathcal{O}(h),
\end{equation}
et d'apr\`es (\ref{pts.6}),
\begin{equation*}
i(1+\mathcal{O}(h))u_{-1}-i(1+\mathcal{O}(h))u_{0}+i(1+\mathcal{O}%
(h))u_{1}=0,
\end{equation*}
d'o\`u
\begin{equation}  \label{qbss.15}
u_0=(1+\mathcal{O}(h))u_1+(1+\mathcal{O}(h))u_{-1}.
\end{equation}

Pour $z\in ]\alpha _{0}^{0},\beta _{0}^{0}[$, $E=E_{0}$, $\varepsilon =0$,
comparons (cf la remarque \ref{pts1})
\begin{equation}
\begin{split}
& u_{-1}=(1+\mathcal{O}(h))(i\varphi _{-1,\alpha _{0}}^{\prime })^{-\frac{1}{%
2}}e^{i\varphi _{-1,\alpha _{0}}/h} \\
& \text{et} \\
& u_{1}=(1+\mathcal{O}(h))(i\varphi _{1,\alpha _{0}}^{\prime })^{-\frac{1}{2}%
}e^{i\varphi _{1,\alpha _{0}}/h}.
\end{split}
\label{qbss.16}
\end{equation}%
D'apr\`{e}s (\ref{qbss.9}), nous avons $\varphi _{-1,\alpha _{0}}^{\prime
}>0 $ (aussi $\varphi _{1,\alpha _{0}}=-\varphi _{-1,\alpha _{0}}$) et donc,
\begin{equation*}
\mathrm{arg\,}(i\varphi _{-1,\alpha _{0}}^{\prime })^{1/2}\in \{\pi /4,-3\pi
/4\},\ \mathrm{arg\,}(i\varphi _{1,\alpha _{0}}^{\prime })^{1/2}\in \{-\pi
/4,+3\pi /4\}.
\end{equation*}%
Comme $\nu _{-1,1}=1$ nous avons aussi
\begin{equation}
(i\varphi _{-1,\alpha _{0}}^{\prime })^{1/2}=i(i\varphi _{1,\alpha
_{0}}^{\prime })^{1/2}.  \label{qbss.17}
\end{equation}%
Les seules possibilit\'{e}s sont alors
\begin{equation*}
(\mathrm{arg\,}(i\varphi _{-1,\alpha _{0}}^{\prime })^{1/2},\mathrm{arg\,}%
(i\varphi _{1,\alpha _{0}}^{\prime })^{1/2})=(\frac{\pi }{4},-\frac{\pi }{4})%
\text{ ou }(-\frac{3\pi }{4},\frac{3\pi }{4}).
\end{equation*}%
\`{A} des facteurs $1+\mathcal{O}(h)$ pr\`{e}s, on voit alors de (\ref%
{qbss.16}) que $u_{-1}=\overline{u}_{1}$ et (\ref{qbss.15}) est bien en
accord avec le fait que $u_{0}$ est une solution r\'{e}elle de l'\'{e}%
quation de Schr\"{o}dinger.

Voici une fa\c con plus directe de d\'eterminer $\mathrm{arg\,}(i\varphi
^{\prime}_{\mp 1,\alpha _0})^{1/2}$: Pour $z<\alpha _0$ proche de $\alpha _0$%
, nous avons $(i\varphi ^{\prime}_{0,\alpha _0})^{1/2}=(V_0(z)-E_0)^{1/4}$,
la branche principale positive. Tournons maintenant dans le sens n\'egatif
vers $\Sigma _{-1}$. Dans la r\'egion de transition entre les deux secteurs $%
\Sigma _0$ et $\Sigma _{-1}$ nous avons, puisque $\nu _{-1,0}=1$, $(i\varphi
^{\prime}_{-1,\alpha _0})^{1/2}=i(V_0(z)-E_0)^{1/4}$ avec la m\^eme branche
de la racine quatri\`eme, o\`u on met une coupure le long de $]\alpha
_0^0,\beta _0^0[$. Quand on arrive \`a $]\alpha_0^0,\beta _0^0[$, on a donc
\begin{equation*}
\mathrm{arg\,}(i\varphi ^{\prime}_{1,\alpha _0})^{1/2}=\frac{\pi }{2}+%
\mathrm{arg\,}(V_0(z)-E_0)^{1/4}=\frac{\pi }{2}-\frac{\pi }{4}=\frac{\pi }{4}%
:
\end{equation*}
\begin{equation}  \label{qbss.18}
\mathrm{arg\,}(i\varphi ^{\prime}_{\mp 1,\alpha _0})^{1/2}=\pm \frac{\pi }{4}%
,\text{ pour }E=E_0,\ \varepsilon =0,\ \alpha _0^0<z<\beta _0^0.
\end{equation}

Regardons maintenant les $v_{j}$, qui en analogie avec (\ref{qbss.3})
prennent la forme
\begin{equation}
v_{j}=\frac{1+\mathcal{O}(h)}{(-i\varphi _{j,\beta _{0}}^{\prime })^{1/2}}%
e^{i\varphi _{j,\beta _{0}}/h}\text{ dans }\mathrm{int\,}(S_{j}).
\label{qbss.19}
\end{equation}%
Nous avons introduit le signe $-$ dans la racine carr\'{e} car nous avons
maintenant $-i\varphi _{0,\beta _{0}}>0$ quand $E$ est r\'{e}el, $%
\varepsilon =0$ et $z>\beta _{0}$. (Ceci se comprend encore mieux si on
travaille avec la variable $-z$ \`{a} la place de $z$). Ici (cf. (\ref%
{qbss.8}), (\ref{qbss.9})) nous avons pour $j=\pm 1$,
\begin{equation}
\varphi _{j,\alpha _{0}}=\varphi _{j,\beta _{0}}+\varphi _{j,\alpha
_{0}}(\beta _{0})\hbox{ pour }z\in S_{j}.  \label{qbss.20}
\end{equation}%
Comme pour les $u_{j}$ il faut discuter le choix de la racine carr\'{e} de $%
(-i\varphi _{j,\beta _{0}})^{1/2}$. Pour cela on choisit la branche de $%
(-\varphi _{0,\beta _{0}}^{\prime })^{1/2}$ qui est positive $=(V-E)^{1/4}$
quand $E$ est r\'{e}el, $\varepsilon =0$ et $z>\beta _{0}$. On d\'{e}finit
ensuite les nombres $\mu _{j,k}$ par
\begin{equation}
(-i\varphi _{j,\beta _{0}}^{\prime })^{1/2}=i^{\mu _{j,k}}(-i\varphi
_{k,\beta _{0}}^{\prime })^{1/2}\hbox{ dans }\mathrm{int\,}(S_{j}\cup
S_{k}),\ j\neq k.  \label{qbss.21}
\end{equation}%
On fait le m\^{e}me choix des $\mu _{j,k}$ que des $\nu _{j,k}$, en
progressant de $S_{0}$ \`{a} la place de $\Sigma _{0}$ dans le sens positif:
\begin{equation}
\mu _{j,k}=\nu _{-j,-k}.  \label{qbss.22}
\end{equation}%
Ainsi,
\begin{equation}
\mu _{-1,0}=-1,\ \mu _{0,1}=-1,\ \mu _{1,-1}=1,\quad \mu _{k,j}=-\mu _{j,k}.
\label{qbss.23}
\end{equation}

Alors en analogie avec (\ref{qbss.15}) et (\ref{qbss.18}) nous avons
\begin{equation}  \label{qbss.24}
v_0=(1+\mathcal{O}(h))v_1+(1+\mathcal{O}(h))v_{-1},
\end{equation}
\begin{equation}  \label{qbss.25}
\mathrm{arg\,}(-i\varphi ^{\prime}_{\pm 1,\beta _0})^{1/2}=\pm \frac{\pi }{4}%
, \hbox{ dans }]\alpha _0,\beta _0[, \hbox{ quand }E\in \mathrm{vois\,}(E_0,%
\mathbb{R}),\ \varepsilon =0.
\end{equation}

Quitte \`{a} modifier $v_{j}$ et $u_{j}$ par des facteurs constants $1+%
\mathcal{O}(h)$, on peut supposer que
\begin{equation}
u_{0}=u_{1}+u_{-1},\ v_{0}=v_{1}+v_{-1}.  \label{qbss.26}
\end{equation}%
remarquons aussi que quand $E$ est r\'{e}el et $\varepsilon =0$, alors $%
u_{-j}=(1+\mathcal{O}(h))\overline{u}_{j}$ et de m\^{e}me pour $v_{j}$.
(Nous allons revenir \`{a} cette propri\'{e}t\'{e} et montrer que les pr\'{e}%
facteurs peuvent \^{e}tre \'{e}limin\'{e}s).

Se rappelant que $u_j$ et $v_j$ sont colin\'eaires pour $j=\pm 1$, en
utilisant (\ref{qbss.26}) on obtient,
\begin{equation}  \label{qbss.27}
W(u_0,v_0)=W(u_{-1},v_1)+W(u_1,v_{-1}).
\end{equation}
Ici, on remarque que
\begin{equation*}
\begin{split}
\varphi _{-1,\alpha _0}(z)+\varphi _{1,\beta _0}(z)&=\int_{\alpha _0}^z
(E-V_{\varepsilon }(t))^{1/2}dt-\int_{\beta _0}^z (E-V_{\varepsilon
}(t))^{1/2}dt \\
&=\int_{\alpha _0}^{\beta _0}(E-V_{\varepsilon }(t))^{1/2}dt=:\frac{1}{2}%
I(E,\varepsilon ),
\end{split}%
\end{equation*}
o\`u la derni\`ere \'egalit\'e d\'efinit l'action $I(E,\varepsilon )$ et
o\`u on choisit la branche de la racine carr\'e qui est positif quand $E$
est r\'eel, $\varepsilon =0$ et $\alpha _0<t<\beta _0$. De m\^eme,
\begin{equation*}
\varphi _{1,\alpha _0}(z)+\varphi _{-1,\beta _0}(z)=-\int_{\alpha _0}^{\beta
_0}(E-V_{\varepsilon }(t))^{1/2}dt=:\frac{1}{2}I(E,\varepsilon )=-\frac{1}{2}%
I(E,\varepsilon ),
\end{equation*}

Nous pouvons voir $(E-V_{\varepsilon })^{\frac{1}{2}}$, comme une fonction
holomorphe sur $U\backslash \lbrack \alpha _0,\beta _0]$, soit $\gamma $ un
contour ferm\'{e} autour de $[\alpha _0,\beta _0]$ orient\'e dans le sens n%
\'{e}gatif. Alors
\begin{equation}
I(E):=2\int_{\alpha _{0}}^{\beta _{0}}(E-V_{\varepsilon })^{\frac{1}{2}%
}dz=\int_{\gamma }(E-V_{\varepsilon })^{\frac{1}{2}}dz.  \label{qbss.28}
\end{equation}

Quand $E\in \mathbb{R}$, $\varepsilon =0$ nous avons aussi
\begin{equation}
I(E)=\int_{p^{-1}(E)}\xi dx=\mathrm{vol}_{\mathbb{R}\times \mathbb{R}%
}\,p^{-1}(]-\infty ,E[),  \label{qbss.29}
\end{equation}
o\`u $p(x,\xi )=p_{\varepsilon }(x,\xi )=\xi ^2+V(x)$ est le symbole
semi-classique de $P=P_{\varepsilon }$ et o\`u la courbe r\'eelle est
orient\'ee dans la direction du champ hamiltonien $H_p=p^{\prime}_{\xi }
\partial _x-p^{\prime}_x \partial _{\xi } $. Rappelons aussi que pour $E\in
\mathbb{R}$, $\varepsilon =0$,
\begin{equation}
\partial_E I(E)=T(E)>0,  \label{qbss.30}
\end{equation}
est la p\'{e}riode primitive pour le flot de $H_p$ dans la courbe d'\'{e}%
nergie r\'{e}elle $p^{-1}(E) $.

Revenons au calcul de notre Wronskien. On obtient
\begin{equation*}
W(u_{-1},v_{1})=\frac{2i\varphi _{-1,\alpha _{0}}^{\prime }(1+\mathcal{O}(h))%
}{(i\varphi _{-1,\alpha _{0}}^{\prime })^{\frac{1}{2}}(-i\varphi _{1,\beta
_{0}}^{\prime })^{\frac{1}{2}}}e^{\frac{i}{2h}I(E,\varepsilon )}.
\end{equation*}%
Ici,
\begin{equation}
\frac{2i\varphi _{-1,\alpha _{0}}^{\prime }}{(i\varphi _{-1,\alpha
_{0}}^{\prime })^{\frac{1}{2}}(-i\varphi _{1,\beta _{0}}^{\prime })^{\frac{1%
}{2}}}=\pm 1  \label{qbss.31}
\end{equation}%
puisque $-i\varphi _{1,\beta _{0}}^{\prime }=i\varphi _{-1,\alpha
_{0}}^{\prime }$ et pour d\'{e}terminer le signe, on peut se placer sur $%
]\alpha _{0},\beta _{0}[$ en supposant que $E\in \mathbb{R}$, $\varepsilon
=0 $. Alors $i\varphi _{-1,\alpha _{0}}^{\prime }=i(E-V(z))^{1/2}$ est
d'argument $\pi /2$ et nous savons que
\begin{equation*}
\mathrm{arg\,}(i\varphi _{-1,\alpha _{0}}^{\prime })^{\frac{1}{2}}=\frac{\pi
}{4}=\mathrm{arg\,}(-i\varphi _{1,\beta _{0}}^{\prime })^{\frac{1}{2}}.
\end{equation*}%
On a donc $+1$ dans (\ref{qbss.31}):
\begin{equation}
W(u_{-1},v_{1})=2(1+\mathcal{O}(h))e^{\frac{i}{2h}I(E,\varepsilon )}.
\label{qbss.32}
\end{equation}%
De la m\^{e}me fa\c{c}on on obtient
\begin{equation}
W(u_{1},v_{-1})=2(1+\mathcal{O}(h))e^{-\frac{i}{2h}I(E,\varepsilon )}.
\label{qbss.33}
\end{equation}

(\ref{qbss.27}) donne alors
\begin{equation}  \label{qbss.34}
W(u_0,v_0)=2((1+\mathcal{O}(h))e^{\frac{i}{2h}I(E,\varepsilon )}+(1+\mathcal{%
O}(h))e^{-\frac{i}{2h}I(E,\varepsilon )}).
\end{equation}
Ici les facteurs $1+\mathcal{O}(h)$ d\'ependent holomorphiquement de $%
(E,\varepsilon )\in \mathrm{vois\,}((E_0,0),\mathbb{C}^2)$ et ont des
d\'eveloppements asymptotiques en puissances de $h$ dans $\mathrm{Hol\,}(%
\mathrm{vois\,}((E_0,0),\mathbb{C}^2))$. On peut \'ecrire
\begin{equation}  \label{qbss.35}
W(u_0,v_0)=2(1+\mathcal{O}(h))e^{-\frac{i}{2h}I(E,\varepsilon )}\left(e^{%
\frac{i}{h}(I(E,\varepsilon )+h^2r(E,\varepsilon ;h))}+1 \right),
\end{equation}
o\`u
\begin{equation}  \label{qbss.36}
r(E,\varepsilon ;h)\sim r_0(E,\varepsilon )+hr_1(E,\varepsilon )+...%
\hbox{ dans
}\mathrm{Hol\,}(\mathrm{vois\,}(E_0,0),\mathbb{C}^2).
\end{equation}
Les z\'eros de $W(u_0,v_0)$ sont donc donn\'es par la condition de
quantification de Bohr-Sommerfeld,
\begin{equation}  \label{qbss.37}
I(E,\varepsilon )+h^2r(E,\varepsilon ;h)=(k+\frac{1}{2})2\pi h,\ k\in
\mathbb{Z}.
\end{equation}

Rappelons maintenant que
\begin{equation*}
\frac{d}{dE}I(E,\varepsilon )=T(E,\varepsilon )\ne 0
\end{equation*}
est la p\'eriode primitive du champ hamiltonien $H_p=p^{\prime}_{\xi }
\partial _x-p^{\prime}_x\partial _{\xi } $ sur la courbe d'\'energie
complexe $p=E$ restreinte \`a un petit voisinage de la courbe r\'eelle $%
p_{\varepsilon =0}(x,\xi )=E_0$. Alors, pour $h$ petit, on peut appliquer le
th\'eor\`eme des fonctions implicites dans sa version holomorphe pour
conclure que l'application
\begin{equation}  \label{qbss.38}
\mathrm{vois\,}(E_0,\mathbb{C})\ni E \mapsto I(E,\varepsilon
;h):=I(E,\varepsilon )+h^2r(E,\varepsilon ;h)\in \mathrm{vois\,}(I(E_0,0),%
\mathbb{C})
\end{equation}
est bijective avec l'inverse
\begin{equation*}
w\mapsto I^{-1}(w,\varepsilon ;h)
\end{equation*}
tel que
\begin{equation}  \label{qbss.39}
I^{-1}(w,\varepsilon ;h)=I^{-1}(w,\varepsilon )+h^2(I^{-1})_2(w,\varepsilon
)+h^3(I^{-1})_3(w,\varepsilon )+...
\end{equation}
dans l'espace des fonctions holomorphes en $(w,\varepsilon )$ dans un
voisinage de $(I(E_0,0),0)$. Les valeurs propres de $P_{\varepsilon }$
(c.\`a.d. les z\'eros de $W(u_0,v_0)$) dans un voisinage de $E_0$ sont alors
donn\'ees par
\begin{equation}  \label{qbss.40}
E_k(\varepsilon ;h)=I^{-1}((k+\frac{1}{2})2\pi h,\varepsilon ;h)
\end{equation}
pour $k\in \mathbb{Z}$ tels que $(k+\frac{1}{2})2\pi h$ appartient \`a un
voisinage de $I(E_0)$.

Quand $\varepsilon =0$ alors $P_{\varepsilon }$ est autoadjoint et les
valeurs propres sont r\'eelles. On peut en d\'eduire que les termes dans les
d\'eveloppements asymptotiques en puissances de $h$ de $I$ et de $I^{-1}$
sont r\'eels pour $\varepsilon =0$.

\section{Le cas d'un seul puits $\mathcal{PT}$-sym\'etrique}

On fait les hypoth\`{e}ses de simple puits de la section (\ref{qbss}).
Supposons aussi que $P_{\varepsilon }$ soit $\mathcal{PT}$-sym\'{e}trique
quand $\varepsilon >0$ est r\'{e}el (et donc aussi que $B=-A$):
\begin{equation}
\lbrack \mathcal{P}\mathcal{T},P_{\varepsilon }]=0,  \label{sppt.1}
\end{equation}%
o\`{u} $\mathcal{P}f(x)=f(-x)$, $\mathcal{T}f(x)=\overline{f(\overline{x})}$%
. Pour $E$ complexe, on a
\begin{equation}
\mathcal{P}\mathcal{T}(P_{\varepsilon }-E)=(P_{\varepsilon }-\overline{E})%
\mathcal{P}\mathcal{T}.  \label{sppt.2}
\end{equation}%
Soient $u_{0}(x,\varepsilon ,E)$ et $v_{0}(x,\varepsilon ,E)$ des solutions
de $(P_{\varepsilon }-E)u=0$ comme dans la section pr\'{e}c\'{e}dente.
Remarquons que gr\^{a}ce \`{a} (\ref{sppt.2}) on peut choisir $v_{0}$ de la
forme
\begin{equation}
v_{0}(x,\varepsilon ,E)=\mathcal{P}\mathcal{T}u_{0}(x,\varepsilon ,\overline{%
E})=\overline{u_{0}(-\overline{x},\varepsilon ,\overline{E})}.
\label{sppt.3}
\end{equation}%
Comme dans la section pr\'{e}c\'{e}dente on cherche les valeurs propres pr%
\`{e}s de $E_{0}$ comme les z\'{e}ros de la fonction $W(E):=W(u_{0},v_{0})$.
On trouve
\begin{equation*}
\begin{split}
W(E)& =h\partial _{x}u_{0}(0,\varepsilon ,E)v_{0}(0,\varepsilon
,E)-u_{0}(0,\varepsilon ,E)h\partial _{x}v_{0}(0,\varepsilon ,E) \\
& =h\partial _{x}u_{0}(0,\varepsilon ,E)\overline{u_{0}(0,\varepsilon ,%
\overline{E})}+u_{0}(0,\varepsilon ,E)\overline{h\partial
_{x}u_{0}(0,\varepsilon ,\overline{E})},
\end{split}%
\end{equation*}%
et on voit que
\begin{equation}
\overline{W(\overline{E})}=W(E).  \label{sppt.4}
\end{equation}

D'autre part, nous avons (\ref{qbss.34}):
\begin{equation}  \label{sppt.5}
W(E)=2(a(E,\varepsilon ;h)e^{\frac{i}{2h}I(E,\varepsilon )}+b(E,\varepsilon
;h)e^{-\frac{i}{2h}I(E,\varepsilon )},
\end{equation}
o\`u $a,b=1+\mathcal{O}(h)$ ont des d\'eveloppements asymptotiques en
puissances de $h$ dans $\mathrm{Hol\,}(\mathrm{vois\,}(E_0,0))$.

Pour $E$ et $\varepsilon $ r\'{e}els on sait que $W(E)$ est r\'{e}el et donc
les deux termes dans (\ref{sppt.5}) sont mutuellement conjugu\'{e}s
complexes. Il en r\'{e}sulte d'abord que $I(E,\varepsilon )$ est r\'{e}el
(comme \'{e}nonc\'{e} dans la proposition (\ref{intr3}) et si on veut
ensuite prendre $E$ complexe on a
\begin{equation}
I(\overline{E},\varepsilon )=\overline{I(E,\varepsilon )},  \label{sppt.6}
\end{equation}%
toujours avec $\varepsilon $ r\'{e}el. On voit ensuite que
\begin{equation*}
b(E,\varepsilon ;h)=\overline{a(E,\varepsilon ;h)}
\end{equation*}%
quand $E$ et $\varepsilon $ sont r\'{e}els et donc plus g\'{e}n\'{e}ralement
que
\begin{equation}
b(E,\varepsilon ;h)=\overline{a(\overline{E},\varepsilon ;h)},
\label{sppt.7}
\end{equation}%
pour $E$ complexe, toujours avec $\varepsilon $ r\'{e}el.

Explicitons alors (\ref{qbss.35}):
\begin{equation*}
W(E)=2be^{-\frac{i}{2h}I}(\exp \frac{i}{h}(I(E,\varepsilon )+h^{2}r)+1),
\end{equation*}%
o\`{u}
\begin{equation*}
r=\frac{1}{ih}\ln \frac{a(E,\varepsilon ;h)}{b(E,\varepsilon ;h)}
\end{equation*}%
est r\'{e}el pour $E$ r\'{e}el (toujours avec $\varepsilon $ r\'{e}el). Il
est alors clair que les valeurs propres de $P_{\varepsilon }$ pr\`{e}s de $%
E_{0}$, donn\'{e}es par (\ref{qbss.37}), sont r\'{e}elles. Ceci termine
aussi la preuve du Th\'{e}or\`{e}me (\ref{ThSP}).

\end{document}